\documentclass[11pt]{amsart}

\usepackage{fullpage}
\usepackage{amsmath}
\usepackage{amssymb}
\usepackage{verbatim}
\usepackage[linktocpage]{hyperref}

\newtheorem{df}{Definition}[section]
\newtheorem{thm}[df]{Theorem}

\newtheorem{prop}[df]{Proposition}
\newtheorem{rem}[df]{Remark}

\newtheorem{lem}[df]{Lemma}

\newcommand{\pf}{\textit{Proof.} }

\begin{document}

\title{A Curvature Flow Unifying Symplectic Curvature Flow And Pluriclosed Flow\footnotemark[1]}

\author{Song Dai}
%\footnotemark[1]

\footnotetext[0]{The author is partially supported by China Scholarship Council.}

\address{School of Mathematical Science\\
Peking University\\
No.5 Yiheyuan Road Haidian District\\
Beijing\\
P.R.China 100871}

\email{daisong0620@gmail.com}

\begin{abstract}
Streets and Tian introduced pluriclosed flow and symplectic curvature flow in \cite{ST pluri} and \cite{ST symp}.
Here we construct a curvature flow to unify these two flows.
We show the short time existence of our flow and exhibit an obstruction to long time existence.
\end{abstract}

\maketitle

\tableofcontents

\section{Introduction}
In recent years, Streets and Tian initialized the study of special geometric structures,
such as, generalized K\"{a}hler and symplectic structure,
by using curvature flows they introduced. They include Hermitian curvature flow, pluriclosed flow, almost Hermitian curvature flow and symplectic curvature flow,
\cite{ST pluri}\cite{ST her}\cite{ST symp}. Subsequently, there are several further works along this direction,
see \cite{B}\cite{EFV}\cite{E}\cite{F}\cite{LW}\cite{Po}\cite{Sm}\cite{Streets pgs}\cite{Streets pbg}\cite{ST reg pluri}\cite{ST gkg}\cite{SW}\cite{V}. Especially in \cite{Streets pgs}\cite{Streets pbg}\cite{ST reg pluri}\cite{ST gkg}\cite{SW}, pluriclosed flow is showed to have rich beautiful results and be related to generalized K\"{a}hler geometry.
In this paper, we introduce a curvature flow which unifies symplectic curvature flow and pluriclosed flow.
%For other related results, one may refer [?].

In \cite{ST symp}, Streets and Tian introduced symplectic curvature flow as follows, which preserves almost K\"{a}hler structure,
\setlength\arraycolsep{2pt}
\begin{eqnarray}
\frac{\partial}{\partial t}g&=&-2Ric+\frac{1}{2}B^{1}-B^{2}\nonumber\\
\frac{\partial}{\partial t}J&=&\triangle{J}+\mathcal{N}+\mathcal{R}\label{symp}\\
g(0)&=&g_{0}\nonumber\\
J(0)&=&J_{0},\nonumber
\end{eqnarray}
where $\mathcal{R}$ is curvature term and $B^{1},B^{2},\mathcal{N}$ are all quadratic terms of $DJ$.
We will give the precise definitions of these tensors in Section \ref{secmain}.

In \cite{ST pluri}, they introduced pluriclosed flow as follows, which preserves pluriclosed structure,
\setlength\arraycolsep{2pt}
\begin{eqnarray*}
\frac{\partial}{\partial t}\omega&=&\partial\partial^{*}\omega+\overline\partial\overline\partial^{*}\omega
+\frac{\sqrt{-1}}{2}\partial\overline\partial \text{logdet}g\\
\omega(0)&=&\omega_{0}.
\end{eqnarray*}
Then, in \cite{ST reg pluri}\cite{ST gkg} they observed that after a gauge transformation induced by Lee form $\theta=-Jd^{*}\omega$,
pluriclosed flow is equivalent to the following flow,
\setlength\arraycolsep{2pt}
\begin{eqnarray}
\frac{\partial}{\partial t}g&=&-2Ric+\frac{1}{2}\mathcal{B}\nonumber\\
\frac{\partial}{\partial t}J&=&\triangle{J}+\mathcal{R}+\mathcal{Q}\label{pluri}\\
g(0)&=&g_{0}\nonumber\\
J(0)&=&J_{0}\nonumber,
\end{eqnarray}
where $\mathcal{B}$ and $\mathcal{Q}$ are quadratic terms of $DJ$. We will give the precise definitions of these tensors in Section \ref{secmain}.
Upon this setting, in \cite{ST gkg}, they showed that twisted generalized K\"{a}hler manifold is a natural background to run pluriclosed flow.

In \cite{H}, Hitchin first introduced the notion of generalized complex structure, which unifies symplectic structure and complex structure.
After that Gualtieri discussed generalized complex structure in detail in his thesis \cite{G}.
In \cite{G}, Gualtieri discovered that a pair of compatible almost generalized complex structures $(\mathcal{J}_{1}, \mathcal{J}_{2})$
is equivalent to almost biHermitian data $(g, J_{+}, J_{-}, b)$, where $J_{\pm}$ are almost complex structures,
compatible with $g$, and $b$ is a 2-form. If $\mathcal{J}_{1}, \mathcal{J}_{2}$ are both integrable, i.e. generalized K\"{a}hler,
the integrability condition is equivalent to
\begin{eqnarray*}
N_{J_{+}}=N_{J_{-}}=0,\\
-d^{c}_{+}\omega_{+}=d^{c}_{-}\omega_{-}=db.
\end{eqnarray*}
If we only require $db$ to be a closed 3-form $H$, which is the twisted case, Streets and Tian \cite{ST gkg} showed that the equivalent pluriclosed flow (\ref{pluri}) of $(g,J_{+})$ and $(g,J_{-})$ preserves generalized K\"{a}hler structure.

A symplectic structure $\omega$ gives a generalized complex structure $\mathcal{J}_{\omega}$, and an almost K\"{a}hler structure $(\omega, J)$ gives a compatible pair of almost generalized complex structures ($\mathcal{J}_{\omega}, \mathcal{J}_{J}$), where $\mathcal{J}_{\omega}$ is integrable while $\mathcal{J}_{J}$ is not necessarily. So one may also regard symplectic curvature flow as a curvature flow to deform a compatible pair of almost generalized complex structures ($\mathcal{J}_{1}, \mathcal{J}_{2}$), where $\mathcal{J}_{1}$ is pintegrable. It leads to the problem whether or not there is a curvature flow which unifies the flow in (1) and (2).
The following theorem gives a solution to this problem.
\begin{thm}\label{main}
Let $(M,g_{0},J_{0})$ be an almost Hermitian manifold. Suppose $M$ is compact.
Then there exists a unique family of almost Hermitian structures $(g(t),J(t)),t\in [0,\epsilon)$ on $M$ satisfying the equations:
\setlength\arraycolsep{2pt}
\begin{eqnarray}
\frac{\partial}{\partial t}g&=&-2Ric+Q_{1}\nonumber\\
\frac{\partial}{\partial t}J&=&\triangle{J}+\mathcal{N}+\mathcal{R}+Q_{2}\label{flow}\\
g(0)&=&g_{0}\nonumber\\
J(0)&=&J_{0}.\nonumber
\end{eqnarray}
Here $\mathcal{R}$ and $\mathcal{N}$ are the same as in (\ref{symp}),
$Q_{1}$ and $Q_{2}$ are quadratic terms of $DJ$ (See Section \ref{secmain} for their precise definitions).
This flow preserves the integrability of $J$.
Furthermore, if the initial data is almost K\"{a}hler, this flow coincides with symplectic curvature flow and if the initial data is pluriclosed,
this flow is equivalent to pluriclosed flow. In particular, if the initial data is K\"{a}hler, this flow is K\"{a}hler-Ricci flow.
\end{thm}

Another motivation to unify (\ref{symp}) and (\ref{pluri}) is to try to understand symplectic curvature flow better.
The tremendous success of Perelman's work \cite{Pe} motivates us to consider finding similar tools in symplectic curvature flow
as in Ricci flow. To begin with, we consider whether symplectic curvature flow is a gradient flow as Ricci flow.
It seems difficult to construct such a functional directly. But as shown in \cite{ST reg pluri},
pluriclosed flow is a gradient flow and the functional is similar to the case of Ricci flow.
So maybe our flow could give some hints to discover the desired functional in symplectic curvature flow.

Turning to regularity, we derive the evolution equations, and then obtain the derivative estimates as follows.
\begin{thm}\label{main2}
Let $(M, g(t), J(t))$ be a solution of (\ref{flow}) for $t\in[0,T)$. Suppose $M$ is compact. If there exists a constant $K$, such that
\begin{eqnarray*}
\sup_{\lbrack 0,T)\times M}
\{t|Rm|, t^{\frac{1}{2}}|DJ|\}\leq K,
\end{eqnarray*}
then for $k\geq 0$, there exists a constant $C=C(k,n,K)$, such that
\begin{eqnarray*}
\sup_{[0,T)\times M}\{t^{\frac{k+2}{2}}|D^{k}Rm|, t^{\frac{k}{2}}|D^{k}J|\}\leq C.
\end{eqnarray*}
\end{thm}

Finally, we obtain an obstruction to long time existence.
\begin{thm}\label{main3}
Let $(M, g(t), J(t))$ be a solution of (\ref{flow}) for $t\in[0,T)$, $T$ be the maximal existence time, $T<+\infty$. Suppose $M$ is compact. Then
\begin{eqnarray*}
\sup_{\lbrack 0,T)\times M}
\{|Rm|, |DJ|\}=+\infty.
\end{eqnarray*}
\end{thm}

%\textbf{Outline of the proof:}
We outline the proof now. Some results in this paper can be implied directly from the results in \cite{ST symp}. For the convenience of readers, we give the complete proof here.
%The proof is similar to \cite{ST symp}.

To prove Theorem \ref{main}, we use DeTurck trick. But we notice that the almost complex structure $J$ does not live in a vector space.
So we transform the equation on the space of almost complex structure to its tangent space at $J_{0}$.
We don't assume $(g,J)$ is compatible first, so we do some modifications to ensure the compatibility,
which gives the non-degenerate symbol. Thus we obtain the short time existence of the modified flow.
Then we do some estimates to show the modified flow gives compatible pair $(g,J)$ and coincide with the initial flow.
For uniqueness, it is the same as in Ricci flow. In the case of symplectic setting and pluriclosed setting,
by direct calculation in Section \ref{secmain}, we see that this flow can be reduced to symplectic curvature flow and pluriclosed flow respectively. So by uniqueness, they coincide with our flow. And the similar argument also applies to the integrability of $J$.

To prove Theorem \ref{main2}, the argument is standard. We derive the evolution equations of $D^{k}J$ and $D^{k}Rm$,
then we construct a function involving the terms we want to estimate. Calculating the evolution equation of this function,
and then by maximum principle, we obtain the desired result. To prove Theorem \ref{main3}, the argument is also standard and the same as in Ricci flow.

%\textbf{Plan of the paper:}
We organize the paper as follows.
In Section \ref{pre}, we recall some preliminaries in almost Hermitian geometry and derive the necessary condition of a variation of almost Hermitian pairs.
In Section \ref{secmain}, We define the tensors we will use in this paper.
Then we do some calculations to show our flow satisfies the necessary condition. And also by calculation,
we show the additional tensors will vanish in special cases. In Section \ref{pf1}, we prove Theorem \ref{main}.
In Section \ref{pf23}, we prove Theorem \ref{main2} and Theorem \ref{main3}.
%\begin{df}
\\
\textbf{Acknowledgements:} The author wishes to express his gratitude to his advisor Gang Tian,
for suggesting the author to the problem of constructing new curvature flows preserving generalized complex structure,
encouraging the author all the time and many helpful discussions. The author would also like to thank Jeffrey Streets for his helpful comments and suggestions, especially for pointing out that this flow may also preserve the integrability of $J$. The author would also like to thank CSC and TRAM for supporting the author visiting Princeton University.

\section{Preliminaries}\label{pre}
We fix some conventions first.\\
%\begin{con}
\textbf{Convention:}\\
(\romannumeral1) Let $g$ be a Riemannian structure. We identify $T\in \Gamma(End(TM))$ and $T\in \Gamma(T^{*}M\bigotimes T^{*}M)$ by
\begin{displaymath}
g(T(X),Y)=T(X,Y).
\end{displaymath}
We implicitly use this identification throughout this paper.\\
(\romannumeral2) For the same index, we always mean to take trace with respect to these two positions, i.e.,
to choose an orthonormal basis and take sum.\\
(\romannumeral3) We write $DJ^{*3}$ for $DJ*DJ*DJ$, etc.\\
(\romannumeral4) Sometimes, we write $i$ instead of $e_{i}$ for short.\\
(\romannumeral5) Sometimes, we omit time-parameter $t$ if there is no ambiguity.\\
(\romannumeral6) Let $D$ denote the Levi-Civita connection. And we always use Levi-Civita throughout this paper.
%\end{con}
%\section{Conventions and Definitions}

We come back to the preliminaries.

Let $M$ be a manifold, $J$ be a section of $End(TM)$. We call $J$ an almost complex structure if $J^{2}=-1$. An almost complex structure $J$ is called integrable if $J$ is induced by holomorphic coordinates. By the theorem of Newlander-Nirenberg \cite{NN}, $J$ is integrable if and only if $N=0$, where
\begin{eqnarray*}
N(X,Y)=[JX,JY]-[X,Y]-J[JX,Y]-J[X,JY]
\end{eqnarray*}
is called Nijenhuis tensor.

We call $(g,J)$ an almost Hermitian structure if $g$ is a Riemannian metric, $J$ is an almost complex structure and $(g,J)$ is compatible, where $(g,J)$ is compatible means that
\begin{eqnarray*}
g(JX,JY)=g(X,Y).
\end{eqnarray*}
For almost Hermitian structure $(g,J)$, we define
\begin{eqnarray*}
\omega(X,Y)=g(JX,Y).
\end{eqnarray*}
Moreover, if $J$ is integrable, $(g,J)$ is called a Hermitian structure. If $d\omega=0$, $(g,J)$ is called an almost K\"{a}hler structure. If $J$ is integrable and $d\omega=0$, $(g,J)$ is called a K\"{a}hler structure. If $J$ is integrable and $dd^{c}\omega=0$, $(g,J)$ is called a pluriclosed structure or SKT structure (strong K\"{a}hler with torsion), where
\begin{eqnarray*}
d^{c}\omega(X,Y,Z):=-d\omega(JX,JY,JZ).
\end{eqnarray*}
%\end{df}
\begin{df}
Let $h\in \Gamma(T^{*}M\bigotimes T^{*}M)$. We define
\setlength\arraycolsep{2pt}
\begin{eqnarray*}
h^{sym}(X,Y)&=&\frac{1}{2}(h(X,Y)+h(Y,X)), \\
h^{skew}(X,Y)&=&\frac{1}{2}(h(X,Y)-h(Y,X)).
\end{eqnarray*}
\end{df}
\begin{df}
Let $(g,J)$ be an almost Hermitian structure. Let $h\in \Gamma(T^{*}M\bigotimes T^{*}M)$. We define
\setlength\arraycolsep{2pt}
\begin{eqnarray*}
h^{(1,1)}(X,Y)&=&\frac{1}{2}(h(X,Y)+h(JX,JY)),\\
h^{(0,2)+(2,0)}(X,Y)&=&\frac{1}{2}(h(X,Y)-h(JX,JY)).
\end{eqnarray*}
We call $h$ is $(1,1)$(resp. $(0,2)+(2,0)$) if $h^{(0,2)+(2,0)}=0$(resp. $h^{(1,1)}=0$).
\end{df}

%\section{Basic Material about Almost Hermitian Geometry}
In Lemma \ref{01} and Lemma \ref{03}, we derive the necessary condition of a variation of almost Hermitian pair.
\begin{lem}\label{01}
Let $J_t$ be a family of almost complex structures, $\frac{\partial}{\partial t}J=K$. Then
\begin{eqnarray*}
KJ+JK=0.
\end{eqnarray*}
\end{lem}
\pf
By definition,
\begin{eqnarray*}
0=\frac{\partial}{\partial t}J^{2}=KJ+JK.
\end{eqnarray*}
\qed
\begin{lem}\label{02}
Let $(g,J)$ be an almost Hermitian structure, $K\in \Gamma(End(TM))$. Then
\begin{center}
$KJ+JK=0\Longleftrightarrow K$ is $(0,2)+(2,0)$.
\end{center}
\end{lem}
\pf
By definition,
\begin{eqnarray*}
\langle(KJ+JK)X,Y\rangle=K(JX,Y)-K(X,JY)=2K^{(1,1)}(JX,Y).
\end{eqnarray*}
\qed
\begin{rem}
Similarly, $KJ=JK\Longleftrightarrow K$ is $(1,1)$.
\end{rem}
\begin{lem}\label{03}
Let $J_t$ be a family of almost complex structures, $\frac{\partial}{\partial t}J=K$. Let $g_t$ be a family of Riemannian structures compatible with $J_t$, $\frac{\partial}{\partial t}g=h$. Then
\begin{eqnarray*}
K^{sym}J=h^{(0,2)+(2,0)}.
\end{eqnarray*}
\end{lem}
\pf
By using $KJ+JK=0$,
\setlength\arraycolsep{2pt}
\begin{eqnarray*}
0&=&\frac{\partial}{\partial t}(g(JX,JY)-g(X,Y))\\
&=&h(JX,JY)-h(X,Y)+g(KX,JY)+g(JX,KY)\\
&=&-2h^{(0,2)+(2,0)}(X,Y)+K(JX,Y)+K(Y,JX)\\
&=&-2h^{(0,2)+(2,0)}(X,Y)+2(K^{sym}J)(X,Y).
\end{eqnarray*}
\qed
\begin{lem}\label{04}
Let $(g,J)$ be an almost Hermitian structure. Then $(L_{X}g,\:L_{X}J)$ satisfies the necessary condition of a variation of $(g,J)$, i.e.\\
{\renewcommand\baselinestretch{1.2}\selectfont
(\romannumeral1) $L_{X}g$ is symmetric,\\
(\romannumeral2) $L_{X}J$ is $(0,2)+(2,0)$\\
(\romannumeral3) $(L_{X}J)^{sym}J=(L_{X}g)^{(0,2)+(2,0)}$.
\par}
\end{lem}
\pf
Let $\phi_{t}$ be the 1-parameter transformation groups generated by $X$, $g_{t}=\phi_{t}^{\ast}g$, $J_{t}=\phi_{t}^{\ast}J$, then
\begin{eqnarray*}
\frac{\partial}{\partial t}\Big{|}_{t=0}g_{t}=L_{X}g,\quad
\frac{\partial}{\partial t}\Big{|}_{t=0}J_{t}=L_{X}J.
\end{eqnarray*}
Then Lemma \ref{04} follows from Lemma \ref{01}, \ref{02}, \ref{03}.
\qed

\begin{lem}\label{J}
Let $(g,J)$ be an almost Hermitian structure. Then
\begin{eqnarray*}
\langle(D_{X}J)Y,Z\rangle&=&-\langle(D_{X}J)Z,Y\rangle,\\
(D_{X}J)JY&=&-J(D_{X}J)Y.
\end{eqnarray*}
\end{lem}
\pf
Let $X,Y,Z$ be in a normal coordinate system, then
%\begin{eqnarray*}
%{\renewcommand\baselinestretch{1.2}\selectfont
\begin{eqnarray*}
\langle(D_{X}J)Y,Z\rangle&=&\langle D_{X}(JY),Z\rangle=X\langle JY,Z\rangle=-X\langle Y,JZ\rangle=-\langle(D_{X}J)Z,Y\rangle,\\
%\end{eqnarray*}
%So $D_{X}J$ is skew.
%\begin{eqnarray*}
(D_{X}J)JY&=&D_X(JJY)-JD_{X}(JY)=-J(D_{X}J)Y.
\end{eqnarray*}
%\par}
%\end{eqnarray*}
%So $D_{X}J$ is $(0,2)+(2,0)$ from Lemma \ref{02}.
\qed
\begin{lem}\label{gau}
Let $(g,J)$ be an almost Hermitian structure, then
\begin{eqnarray*}
\langle(D_{JX}J)Y,Z\rangle-\langle J(D_{X}J)Y,Z\rangle&=&\frac{1}{2}(N(X,Y,Z)+N(Z,X,Y)-N(Y,Z,X)),\\
\langle(D_{JX}J)Y,Z\rangle+\langle J(D_{X}J)Y,Z\rangle&=&(d\omega)^{+}(JX,Y,Z)-(d\omega)^{+}(JX,JY,JZ).
\end{eqnarray*}
In particular,
\begin{eqnarray*}
D_{JX}J=JD_{X}J&\Longleftrightarrow& N=0\\
D_{JX}J=-JD_{X}J&\Longleftrightarrow& (d\omega)^{+}=0.
\end{eqnarray*}
\end{lem}
\pf
It is proved in \cite{Ga}.
\qed
%\begin{lem}\label{dj2}
%Let $(g,J)$ be an almost Hermitian structure, then
%\begin{eqnarray*}
%\langle(D_{X}J)Y,Z\rangle=\frac{1}{2}(d\omega(X,Y,Z)-d\omega(X,JY,JZ)+\langle X,N(Y,JZ)\rangle).
%\end{eqnarray*}
%\end{lem}
%\pf
%It is proved in \cite{KN}.
%\qed

\section{Main Calculations}\label{secmain}
First, we define the tensors we use in this paper.
\begin{df}\label{defn}
Let $(M,g,J)$ be an almost Hermitian manifold, $X,Y,Z\in TM$.
\setlength\arraycolsep{1pt}
{\renewcommand\baselinestretch{1.2}\selectfont
\begin{eqnarray*}
%Q_{1}(X,Y)&=&-\langle(D_{i}J)X, (D_{i}J)Y\rangle-\langle(D_{X}J)i, (D_{Y}J)i\rangle^{(1,1)}
%-\frac{1}{2}\langle(D_{X}J)i, (D_{Y}J)i\rangle^{(0,2)+(2,0)}\\
%&&-2\langle(D_{i}J)JX,(D_{Ji}J)Y\rangle-\frac{3}{2}\langle(D_{JX}J)i, (D_{Y}J)Ji\rangle^{sym}\\
%&&+2\langle(D_{(D_{i}J)X}J)Y,i\rangle^{(1,1)}+\langle(D_{(D_{i}J)X}J)Y,i\rangle^{(0,2)+(2,0)}\\
%&&+8\langle(D_{X}J)i,(D_{i}J)Y\rangle^{sym,(1,1)}+2\langle(D_{X}J)i,(D_{i}J)Y\rangle^{sym,(0,2)+(2,0)},\\
%Q_{2}(X,Y)&=&-\langle(D_{J(D_{i}J)i}J)X,Y\rangle-\langle(D_{(D_{i}J)i}J)JX,Y\rangle
%+\frac{1}{2}\langle(D_{JX}J)i,(D_{Y}J)i\rangle^{(0,2)+(2,0)}\\
%&&-2\langle(D_{JX}J)Y,(D_{i}J)i\rangle^{(0,2)+(2,0),skew}-\langle(D_{(D_{i}J)JX}J)Y,i\rangle^{(0,2)+(2,0)}\\
%&&-2\langle(D_{JX}J)i,(D_{i}J)Y\rangle^{(0,2)+(2,0),sym}+\langle(D_{JX}J)i,(D_{i}J)Y\rangle^{(0,2)+(2,0),skew}.\\
%N(X,Y)&=&(D_{JX}J)(Y)-(D_{JY}J)(X)-J(D_{X}J)(Y)+J(D_{Y}J)(X),\\
B^{1}(X,Y)&=&\langle(D_{X}J)i,(D_{Y}J)i\rangle,\\
B^{2}(X,Y)&=&\langle(D_{i}J)X,(D_{i}J)Y\rangle,\\
B^{3}(X,Y)&=&\langle(D_{(D_{i}J)X}J)i,Y\rangle=-\langle(D_{i}J)X,j\rangle\langle(D_{j}J)Y,i\rangle\\
B^{4}(X,Y)&=&\langle(D_{X}J)i,(D_{i}J)Y\rangle,\\
\overline{B^{1}}(X,Y)&=&\langle(D_{X}J)i,(D_{Y}J)Ji\rangle,\\
\overline{B^{2}}(X,Y)&=&\langle(D_{i}J)X,(D_{Ji}J)Y\rangle,\\
%\overline{B^{3}}(X,Y)&=&\langle(D_{(D_{i}J)X}J)Ji,Y\rangle\\
Q_{1}&=&-\frac{1}{2}(B^{1})^{(1,1)}-(B^{3})^{(0,2)+(2,0)}+4(B^{4})^{(1,1),sym}-(\overline{B^{1}}J)^{(1,1)}-\overline{B^{2}}J\\
Q_{2}&=&(B^{3})^{(0,2)+(2,0)}J\\
\mathcal{N}&=&B^{2}J\\
\mathcal{R}(X,Y)&=&Ric(JX,Y)+Ric(X,JY),\\
\mathcal{Q}&=&B^{2}J+B^{3}J\\
%Q(X)&=&-(D_{i}J)(D_{JX}J)i-J(D_{(D_{i}J)X}J)i+(D_{i}J)(D_{Ji}J)X\\
%&&-(D_{J(D_{i}J)i}J)X+J(D_{(D_{i}J)i}J)X+(D_{JX}J)(D_{i}J)i-J(D_{X}J)(D_{i}J)i,\\
H(X,Y,Z)&=&d^{c}\omega(X,Y,Z)=-d\omega(JX,JY,JZ),\\
\mathcal{B}(X,Y)&=&H(X,i,j)H(Y,i,j),\\
\theta^{\sharp}&=&-J(D_{i}J)i,\\
\overline{N}(X,Y)&=&\frac{1}{2}(N((D_{i}J)X,i,Y)+N(Y,(D_{i}J)X,i)-N(i,Y,(D_{i}J)X))\\
&&-\frac{1}{2}(N(i,(D_{X}J)i,Y)+N(Y,i,(D_{X}J)i)-N((D_{X}J)i,Y,i))-(D_{i}J)N(X,i),\\
\mathcal{K}(X)&=&(D_{i}N)(Ji,X)\\
(d\omega)^{+}(X,Y,Z)&=&\frac{1}{4}(3d\omega(X,Y,Z)+d\omega(JX,JY,Z)+d\omega(JX,Y,JZ)+d\omega(X,JY,JZ)).
\end{eqnarray*}}
\end{df}
The lemmas below are preparation for the proof of Theorem \ref{main}.
\begin{lem}\label{compa}
Let $(g,J)$ be an almost Hermitian structure, then
 $(-2Ric+Q_{1},\triangle J+\mathcal{N}+\mathcal{R}+Q_{2})$ satisfies the necessary condition of a variation.
\end{lem}
\pf
First, we show that $(-2Ric,\triangle J+\mathcal{N}+\mathcal{R})$ satisfies the necessary condition. We need to check the following things:\\
(\romannumeral1) $Ric$ is symmetric,\quad
(\romannumeral2) $\triangle{J}+\mathcal{N}$ is $(0,2)+(2,0)$,\quad
(\romannumeral3) $\mathcal{R}$ is $(0,2)+(2,0)$,\quad
(\romannumeral4) $\triangle{J}$ is skew,\quad
(\romannumeral5) $\mathcal{N}$ is skew,\quad
(\romannumeral6) $\mathcal{R}$ is symmetric,\quad
(\romannumeral7) $\mathcal{R}J=-2Ric^{(0,2)+(2,0)}$.\quad\\
By definition, it is easy to see (\romannumeral1), (\romannumeral3), (\romannumeral6), (\romannumeral7).
For (\romannumeral2), we use normal coordinate to calculate the $(1,1)$ part of $\triangle J$, by using Lemma \ref{J},
\setlength\arraycolsep{2pt}
\begin{eqnarray*}
\langle(\triangle J)(JX),JY\rangle&=&\langle(D_{i}DJ)(i,JX),JY\rangle\\
&=&\langle D_{i}((D_{i}J)(JX))-(D_{i}J)(D_{i}(JX)),JY\rangle\\
&=&-\langle D_{i}(J(D_{i}J)X)+(D_{i}J)(D_{i}(JX)),JY\rangle\\
&=&-\langle (D_{i}J)(D_{i}J)X+JD_{i}((D_{i}J)X)+(D_{i}J)(D_{i}J)X,JY\rangle\\
&=&-2\langle (D_{i}J)(JX),(D_{i}J)Y\rangle-\langle (D_{i}D_{i}J)X,Y\rangle\\
&=&-2\mathcal{N}-\langle(\triangle J)X,Y\rangle.
\end{eqnarray*}
So $\mathcal{N}=-(\triangle J)^{(1,1)}$.
For (\romannumeral4), we also use normal coordinate,
\setlength\arraycolsep{2pt}
\begin{eqnarray*}
\langle(\triangle J)X,Y\rangle&=&\langle D_{i}((D_{i}J)X),Y\rangle\\
&=&\partial_{i}\langle(D_{i}J)X,Y\rangle\\
&=&\partial_{i}\langle D_{i}(JX),Y\rangle-\partial_{i}\langle J(D_{i}X),Y\rangle\\
&=&\partial_{i}\partial_{i}\langle JX,Y\rangle-\partial_{i}\langle JX,D_{i}Y\rangle+\partial_{i}\langle D_{i}X,JY\rangle,
\end{eqnarray*}
then we see $\triangle J$ is skew.
And (\romannumeral5) follows from Lemma \ref{J}.\\
Next, we show that $(Q_{1},Q_{2})$ satisfies the necessary condition. In fact, by applying Lemma \ref{J}, we can easily obtain that all terms in $Q_{1}$ are symmetric and all terms in $Q_{2}$ are $(0,2)+(2,0)$. And $Q_{1}^{(0,2)+(2,0)}=(B^{3})^{(0,2)+(2,0)}$. So we finish the proof.
\qed
\begin{lem}\label{s}
Let $(g,J)$ be an almost Hermitian structure. Suppose $d\omega=0$. Then
\begin{eqnarray*}
Q_{1}&=&\frac{1}{2}B^{1}-B^{2},\\
Q_{2}&=&0.
\end{eqnarray*}
\end{lem}
\pf
Since $d\omega=0$, by Lemma \ref{J} and Lemma \ref{gau}, one see $B^{1}$ and $B^{3}$ are $(1,1)$, $\overline{B^{1}}J=B^{1}$, $\overline{B^{2}}J=B^{2}$. Now we prove $B^{4}=\frac{1}{2}B^{1}$. In fact, we notice that
\begin{eqnarray*}
\langle(D_{X}J)Y,Z\rangle+\langle(D_{Y}J)Z,X\rangle+\langle(D_{Z}J)X,Y\rangle=d\omega(X,Y,Z)=0.
\end{eqnarray*}
Thus,
\begin{eqnarray*}
\langle(D_{X}J)i,(D_{i}J)Y\rangle&=&\langle(D_{i}J)Y,j\rangle\langle(D_{X}J)i,j\rangle=-\langle(D_{(D_{X}J)i}J)Y,i\rangle\\
&=&\langle(D_{Y}J)i,(D_{X}J)i\rangle+\langle(D_{i}J)(D_{X}J)i,Y\rangle=B^{1}(X,Y)-\langle(D_{X}J)i,(D_{i}J)Y\rangle.
\end{eqnarray*}
So $\langle(D_{i}J)X,(D_{Y}J)i\rangle=\frac{1}{2}B^{1}(X,Y)$. So we finish the proof.
\qed
\begin{lem}\label{p}
Let $(g,J)$ be an almost Hermitian structure. Suppose $N=0$. Then
\begin{eqnarray*}
Q_{1}&=&\frac{1}{2}\mathcal{B},\\
Q_{2}&=&\mathcal{Q}-\mathcal{N}.
\end{eqnarray*}
\end{lem}
\pf
The proof is by direct calculations based on Lemma \ref{J} and Lemma \ref{gau}.
We notice that $B^{1}$ is $(1,1)$, $B^{3}$ is $(0,2)+(2,0)$. And $\overline{B^{1}}=B^{1}J$, $\overline{B^{2}}=B^{2}J$.
We also have $B^{4}=0$, since
\begin{eqnarray*}
\langle(D_{X}J)i,(D_{i}J)Y\rangle=\langle(D_{X}J)Ji,(D_{Ji}J)Y\rangle=-\langle J(D_{X}J)i,J(D_{i}J)Y\rangle=-\langle(D_{X}J)i,(D_{i}J)Y\rangle.
\end{eqnarray*}
We can calculate $\mathcal{B}$ in terms of $DJ$,
\begin{eqnarray*}
\mathcal{B}(X,Y)=H(X,i,j)H(Y,i,j)=d\omega(JX,Ji,Jj)d\omega(JY,Ji,Jj)=d\omega(JX,i,j)d\omega(JY,i,j).
\end{eqnarray*}
We have
\begin{eqnarray*}
d\omega(JX,i,j)=\langle(D_{JX}J)i,j\rangle+\langle(D_{Ji}J)j,X\rangle+\langle(D_{Jj}J)X,i\rangle.
\end{eqnarray*}
Calculating term by term,
\begin{eqnarray*}
\langle(D_{JX}J)i,j\rangle\langle(D_{JY}J)i,j\rangle&=&\langle(D_{X}J)i,(D_{Y}J)i\rangle=B^{1}(X,Y)\\
\langle(D_{Ji}J)j,X\rangle\langle(D_{Ji}J)j,Y\rangle=\langle(D_{Jj}J)X,i\rangle\langle(D_{Jj}J)Y,i\rangle&=&\langle(D_{i}J)X,(D_{i}J)Y\rangle=B^{2}(X,Y)\\
\langle(D_{JX}J)i,j\rangle\langle(D_{Ji}J)j,Y\rangle=\langle(D_{JX}J)i,j\rangle\langle(D_{Jj}J)Y,i\rangle&=&-\langle(D_{X}J)i,(D_{i}J)Y\rangle=0\\
\langle(D_{JY}J)i,j\rangle\langle(D_{Ji}J)j,X\rangle=\langle(D_{JY}J)i,j\rangle\langle(D_{Jj}J)X,i\rangle&=&-\langle(D_{X}J)i,(D_{i}J)Y\rangle=0\\
\langle(D_{Ji}J)j,X\rangle\langle(D_{Jj}J)Y,i\rangle=\langle(D_{Ji}J)j,Y\rangle\langle(D_{Jj}J)X,i\rangle&=&-\langle(D_{(D_{i}J)X}J)i,Y\rangle=-B^{3}(X,Y)
\end{eqnarray*}
So
\begin{eqnarray*}
\frac{1}{2}\mathcal{B}=\frac{1}{2}B^{1}+B^{2}-B^{3}.
\end{eqnarray*}
Then we obtain the desired result.
\qed
\begin{rem}
In \cite{ST gkg}, $\mathcal{Q}$ is defined as
\begin{eqnarray*}
\mathcal{Q}(X)&=&-(D_{i}J)(D_{JX}J)i-J(D_{(D_{i}J)X}J)i+(D_{i}J)(D_{Ji}J)X\\
&&-(D_{J(D_{i}J)i}J)X+J(D_{(D_{i}J)i}J)X+(D_{JX}J)(D_{i}J)i-J(D_{X}J)(D_{i}J)i.
\end{eqnarray*}
Since $N=0$, it coincide with our definition.
\end{rem}
\begin{lem}\label{Lee}
Let $(g,J)$ be an almost Hermitian structure, then
\begin{eqnarray*}
L_{\theta^{\sharp}}J=\triangle J+\mathcal{Q}+\mathcal{R}+\mathcal{K}+\overline{N}.
\end{eqnarray*}
\end{lem}
\pf
In \cite{ST gkg}, there is a similar formula. But in our case, we don't assume $N=0$.

We use normal coordinate,
\setlength\arraycolsep{2pt}
\begin{eqnarray}
(L_{\theta^{\sharp}}J)X&=&(L_{-J(D_{i}J)i}J)X\nonumber\\
&=&-[J(D_{i}J)i,JX]+J[J(D_{i}J)i,X]\nonumber\\
&=&-D_{J(D_{i}J)i}(JX)+D_{JX}(J(D_{i}J)i)+JD_{J(D_{i}J)i}X-JD_{X}(J(D_{i}J)i)\nonumber\\
&=&-(D_{J(D_{i}J)i}J)X+(D_{JX}J)(D_{i}J)i+JD_{JX}((D_{i}J)i)-J(D_{X}J)(D_{i}J)i+D_{X}((D_{i}J)i)\nonumber\\
&=&-(D_{J(D_{i}J)i}J)X+(D_{JX}J)(D_{i}J)i+J(D_{JX}(D_{i}J))i-J(D_{X}J)(D_{i}J)i+D_{X}(D_{i}J)\nonumber\\
&=&J(D^{2}J)(JX,i,i)+(D^{2}J)(X,i,i)-(D_{J(D_{i}J)i}J)X+(D_{JX}J)(D_{i}J)i-J(D_{X}J)(D_{i}J)i.\label{21}
\end{eqnarray}
By Ricci identity,
\begin{eqnarray}
(D^{2}J)(X,i,i)&=&(D^{2}J)(i,X,i)+(Rm(X,i)J)i\nonumber\\
&=&(D^{2}J)(i,X,i)+Rm(X,i)(Ji)-JRm(X,i)i\nonumber\\
&=&(D^{2}J)(i,X,i)+Rm(X,i)(Ji)-JRic(X).\label{22}
\end{eqnarray}
Similarly,
\begin{eqnarray}
J(D^{2}J)(JX,i,i)&=&J(D^{2}J)(i,JX,i)+JRm(JX,i)(Ji)+Ric(JX).\label{23}
\end{eqnarray}
Notice that
\begin{eqnarray*}
N(X,Y)=(D_{JX}J)Y-(D_{JY}J)X-J(D_{X}J)Y+J(D_{Y}J)X.
\end{eqnarray*}
Hence,
\begin{eqnarray}
J(D^{2}J)(i,JX,i)&=&JD_{i}((D_{JX}J)i)-J(D_{(D_{i}J)X}J)i-JD_{i}(J(D_{X}J)i)+JD_{i}(J(D_{X}J)i)\nonumber\\
&=&JD_{i}((D_{JX}J)i-J(D_{X}J)i)-J(D_{(D_{i}J)X}J)i+J(D_{i}J)(D_{X}J)i-(D^{2}J)(i,X,i)\nonumber\\
&=&JD_{i}((D_{Ji}J)X-J(D_{i}J)X)+JD_{i}(N(X,i))\nonumber\\
&&-J(D_{(D_{i}J)X}J)i+J(D_{i}J)(D_{X}J)i-(D^{2}J)(i,X,i)\nonumber
%&=&JD^{2}J(i,Ji,X)+J(D_{(D_{i}J)i}J)X-J(D_{i}J)(D_{i}J)X+(\triangle J)X+\overline{N}\nonumber\\
%&&-(D_{i}J)N(X,i)-J(D_{(D_{i}J)X}J)i+J(D_{i}J)(D_{X}J)i-(D^{2}J)(i,X,i).\label{24}
\end{eqnarray}
Notice that
\begin{eqnarray*}
JD_{i}(N(X,i))&=&D_{i}(JN(X,i))-(D_{i}J)N(X,i)\\
&=&D_{i}(N(Ji,X))-(D_{i}J)N(X,i)\\
&=&(D_{i}N)(Ji,X)+N((D_{i}J)i,X)-(D_{i}J)N(X,i).
\end{eqnarray*}
So
\begin{eqnarray}
J(D^{2}J)(i,JX,i)&=&J(D^{2}J)(i,Ji,X)+J(D_{(D_{i}J)i}J)X-J(D_{i}J)(D_{i}J)X+(\triangle J)X\nonumber\\
&&+\mathcal{K}(X)+N((D_{i}J)i,X)-(D_{i}J)N(X,i)\nonumber\\
&&-J(D_{(D_{i}J)X}J)i+J(D_{i}J)(D_{X}J)i-(D^{2}J)(i,X,i).\label{24}
\end{eqnarray}
And
\begin{eqnarray}
N((D_{i}J)i,X)&=&(D_{J(D_{i}J)i}J)X-(D_{JX}J)(D_{i}J)i+(D_{(D_{i}J)i}J)JX-(D_{X}J)J(D_{i}J)i\nonumber\\
&=&(D_{J(D_{i}J)i}J)X-(D_{JX}J)(D_{i}J)i-J(D_{(D_{i}J)i}J)X+J(D_{X}J)(D_{i}J)i.
\end{eqnarray}
By resorting Lemma \ref{gau}, we obtain
\begin{eqnarray}
&&\langle-J(D_{(D_{i}J)X}J)i-(D_{(D_{i}J)JX}J)i,Y\rangle\nonumber\\
&=&\langle-J(D_{(D_{i}J)X}J)i+(D_{J(D_{i}J)X}J)i,Y\rangle\nonumber\\
&=&\frac{1}{2}(N((D_{i}J)X,i,Y)+N(Y,(D_{i}J)X,i)-N(i,Y,(D_{i}J)X)),
\end{eqnarray}
and
\begin{eqnarray}
&&\langle J(D_{i}J)(D_{X}J)i,Y\rangle\nonumber\\
&=&\langle J(D_{i}J)(D_{X}J)i-(D_{Ji}J)(D_{X}J)i,Y\rangle\nonumber\\
&=&-\frac{1}{2}(N(i,(D_{X}J)i,Y)+N(Y.i,(D_{X}J)i)-N((D_{X}J)i,Y,i)).
\end{eqnarray}
Then by Ricci identity again,
\begin{eqnarray}
JD^{2}J(i,Ji,X)&=&\frac{1}{2}(JD^{2}J(i,Ji,X)-JD^{2}J(Ji,i,X))\nonumber\\
&=&\frac{1}{2}J(Rm(i,Ji)J)X\nonumber\\
&=&\frac{1}{2}(JRm(i,Ji)(JX)+Rm(i,Ji)X).\label{26}
\end{eqnarray}
By Bianchi identity,
\begin{eqnarray*}
Rm(i,Ji)(JX)+Rm(Ji,JX)i+Rm(JX,i)(Ji)=0.
\end{eqnarray*}
Notice that
\begin{eqnarray*}
Rm(Ji,JX)i=Rm(JX,i)(Ji).
\end{eqnarray*}
Thus
\begin{eqnarray}
JRm(i,Ji)(JX)&=&-2JRm(JX,i)(Ji),\\\label{27}
Rm(i,Ji)(X)&=&-2Rm(X,i)(Ji).\label{28}
\end{eqnarray}
Putting (\ref{21})$\sim$(\ref{28}) together, we obtain the desired result.
\qed

\section{Proof of Theorem \ref{main}}\label{pf1}
\textit{Proof of Theorem \ref{main}. } The argument is the same as in \cite{ST symp}.
We use DeTurck trick to prove short time existence and uniqueness.

We consider the following equations,
\setlength\arraycolsep{2pt}
\begin{eqnarray}
\frac{\partial}{\partial t}g&=&-2Ric+Q_{1}+L_{X}g\triangleq \mathcal{D}_{1}(g,J)\nonumber\\
\frac{\partial}{\partial t}J&=&\triangle{J}+\mathcal{N}+\mathcal{R}+Q_{2}+L_{X}J\triangleq \mathcal{D}_{2}(g,J)\label{17}\\
g(0)&=&g_{0}\nonumber\\
J(0)&=&J_{0}\nonumber
\end{eqnarray}
where $X=tr_{g}(\Gamma-\overline{\Gamma})$,  $\overline{\Gamma}$ is the Christoffel symbol of a fixed metric $\overline{g}$.

Then, in order to use the PDE theory in Banach space, we consider the tangent space at $J_{0}$. Denote $T\mathcal{J}_{J}$ the tangent space at $J$, i.e.
\begin{displaymath}
T\mathcal{J}_{J}=\left\{ E\in End(TM)| EJ+JE=0 \right\}.
\end{displaymath}
Then, in a neighborhood $U$ of $J_{0}$, we can identify $J$ and $E$ by using
\begin{displaymath}
\pi:T\mathcal{J}_{J_{0}}\supset U^{\prime} \rightarrow U, ~~~\pi E=-J_{0}e^{J_{0}E}
\end{displaymath}
and note that $D\pi| _{0}=Id$.

Notice that we don't assume that $(g,J)$ is compatible. So we need to do some modifications.
For convenience, we write $g^{J}$(resp. $g^{-J}$) instead of $g^{(1,1)}$(resp. $g^{(0,2)+(2,0)}$) and do similar things for other tensors.
Note that $g^{J}$ is compatible with $J$. We consider the following equations,
\setlength\arraycolsep{2pt}
\begin{eqnarray}
\frac{\partial}{\partial t}g&=&\mathcal{D}_{1}(g^{\pi E},\pi E)+\triangle_{g_{0}}(g^{-\pi E})\triangleq \tilde{\mathcal{D}}_{1}(g,E)\nonumber\\
\frac{\partial}{\partial t}E&=&(D\pi|_{\pi E})^{-1}\mathcal{D}_{2}(g^{\pi E},\pi E)\triangleq\tilde{\mathcal{D}}_{2}(g,E)\label{16}\\
g(0)&=&g_{0}\nonumber\\
E(0)&=&0.\nonumber
\end{eqnarray}
Note that $\tilde{\mathcal{D}}_{1}$ is symmetric, and $\tilde{\mathcal{D}}_{2}$ is well defined since $\triangle{J}+\mathcal{N}+\mathcal{R}+Q_{2}+L_{X}J$ is $(0,2)+(2,0)$ for the pair $(g^{J},J)$.
So $\tilde{\mathcal{D}}_{1}\bigoplus\tilde{\mathcal{D}}_{2}$ gives an operator from $\Gamma((T^{*}M\bigotimes^{sym}T^{*}M)\bigoplus T\mathcal{J}_{J_{0}})$ to itself.

Now, we calculate the symbol of $\tilde{\mathcal{D}}_{1}\bigoplus\tilde{\mathcal{D}}_{2}$ at $(g_{0}, 0)$ to show the short time existence of the modified flow. First, we calculate the variation of $\tilde{\mathcal{D}}_{1}$ along the direction of $(h,0)$, where $h=\delta g$. Since $\delta E=0$, $\pi E=\pi 0=J_{0}$ is fixed. And note that $\delta(g^{J_{0}})=h^{J_{0}}$, $g_{0}^{J_{0}}=g_{0}$. Therefore
\begin{displaymath}
\mathcal{L}_{(g_{0},0)}(\mathcal{D}_{1}(g^{\pi E},\pi E))(h,0)=\mathcal{L}_{g_{0}^{J_{0}}}(\mathcal{D}_{1}(g,J_{0}))(h^{J_{0}})
=\mathcal{L}_{g_{0}}(\mathcal{D}_{1}(g,J_{0}))(h^{J_{0}})
\end{displaymath}
where $\mathcal{L}_{(g_{0},0)}$ denote the linearization operator at $(g_{0},0)$.

Noting that only $-2Ric$ and $L_{X}g$ involve second order term and from the standard calculations in Ricci flow \cite{CK}, we have
\begin{displaymath}
\mathcal{L}_{g_{0}}(\mathcal{D}_{1}(g,J_{0}))(h^{J_{0}})=\triangle_{g_{0}}(h^{J_{0}})+\mathcal{O}(\partial h).
\end{displaymath}
And
\begin{displaymath}
\mathcal{L}_{(g_{0},0)}(\triangle_{g_{0}}(g^{-J}))(h,0)=\triangle_{g_{0}}(h^{-J_{0}}).
\end{displaymath}
Let $\sigma$ denote the symbol of a linear differential operator. Thus we obtain
\begin{displaymath}
\sigma(\mathcal{L}_{(g_{0},0)}\tilde{\mathcal{D}}_{1})(h,0)(x,\xi)=|\xi|^{2}h, \quad \text{where } \xi\in T^{*}_{x}M.
\end{displaymath}
Then we calculate the variation of $\tilde{\mathcal{D}}_{1}$ along the direction of $(0,K)$, where $K=\delta E$. Since $D\pi|_{0}=Id$, we have
\begin{displaymath}
\delta(\tilde{\mathcal{D}}_{1}(g,E))(0,K)=\delta(\mathcal{D}_{1}(g^{J},J))(0, \delta J).
\end{displaymath}
We identify $\delta J$ and $K$ below.

From the calculations above, we see
\begin{displaymath}
(-2Ric(g^{J})+L_{X(g^{J})}(g^{J}))_{ij}=(g^{J})^{pq}\partial_{p}\partial_{q}(g^{J})_{ij}+\mathcal{O}(\partial g, \partial J).
\end{displaymath}
So
\begin{displaymath}
\mathcal{L}_{(g_{0},0)}(\mathcal{D}_{1}(g^{\pi E},\pi E))(0,K)=\frac{\partial}{\partial t}\Big{|}_{t=0}(g_{0})^{pq}\partial_{p}\partial_{q}(g_{0}^{J_{t}})_{ij}+\mathcal{O}(\partial K).
\end{displaymath}
It is easy to see
\begin{displaymath}
\mathcal{L}_{(g_{0},0)}(\triangle_{g_{0}}(g^{-J}))(0,K)=\frac{\partial}{\partial t}\Big{|}_{t=0}(g_{0})^{pq}\partial_{p}\partial_{q}(g_{0}^{-J_{t}})_{ij}+\mathcal{O}(\partial K).
\end{displaymath}
Thus we obtain
%\begin{displaymath}
%\mathcal{L}_{(g_{0},0)}(\mathcal{D}_{1}(g^{\pi E},\pi E))(0,K)=\mathcal{O}(\partial K).
%\end{displaymath}
\begin{displaymath}
\sigma(\mathcal{L}_{(g_{0},0)}\tilde{\mathcal{D}}_{1})(0,K)(x,\xi)=0, \quad \text{where }\xi\in T^{*}_{x}M.
\end{displaymath}

Next, we calculate the variation of $\tilde{\mathcal{D}}_{2}$ along the direction of $(\delta g,\delta E)=(h,K)$. We have
\begin{displaymath}
\delta(\tilde{\mathcal{D}}_{2}(g,E))(h,K)=\delta(\mathcal{D}_{2}(g^{J},J))(\delta g, \delta J).
\end{displaymath}
%The variation vector on the righthand side is still $(\delta g, \delta J)=(h,K)$.
Noting that in the expression of $\mathcal{D}_{2}$, only $\triangle J$, $L_{X}J$, and $\mathcal{R}$ involve second order term, so we only need to calculate these three terms. We calculate them for the pair $(g,J)$ first.
\setlength\arraycolsep{2pt}
\begin{eqnarray*}
(\triangle J)(e_{k})&=&g^{ij}D^{2}J(e_{i},e_{j},e_{k})\\
&=&g^{ij}D_{i}((D_{j}J)e_{k})+\mathcal{O}(\partial g,\partial J)\\
&=&g^{ij}D_{i}(D_{j}(Je_{k})-JD_{j}e_{k})+\mathcal{O}(\partial g,\partial J)\\
&=&g^{ij}D_{i}(D_{j}(J_{k}^{l}e_{l})-J(\Gamma_{jk}^{p}e_{p}))+\mathcal{O}(\partial g,\partial J)\\
&=&g^{ij}(D_{i}(\partial_{j}J_{k}^{l}e_{l})+D_{i}(J_{k}^{p}\Gamma_{jp}^{l}e_{l})
-D_{i}(\Gamma_{jk}^{p}J_{p}^{l}e_{l}))+\mathcal{O}(\partial g,\partial J)\\
&=&g^{ij}(\partial_{i}\partial_{j}J_{k}^{l}+J_{k}^{p}\partial_{i}\Gamma_{jp}^{l}
-J_{p}^{l}\partial_{i}\Gamma_{jk}^{p})e_{l}+\mathcal{O}(\partial g,\partial J),\\
\\
(L_{X}J)(e_{k})&=&[X,Je_{k}]-J[X,e_{k}]\\
&=&[X^{p}e_{p},J_{k}^{l}e_{l}]-J[X^{p}e_{p},e_{k}]\\
&=&(X^{p}\partial_{p}J_{k}^{l}-J_{k}^{p}\partial_{p}X^{l}+J_{p}^{l}\partial_{k}X^{p})e_{l}\\
&=&g^{ij}(J_{p}^{l}\partial_{k}\Gamma_{ij}^{p}-J_{k}^{p}\partial_{p}\Gamma_{ij}^{l})e_{l}+\mathcal{O}(\partial g,\partial J),\\
\\
\mathcal{R}(e_{k})&=&(J_{k}^{p}Ric_{p}^{l}-J_{p}^{l}Ric_{k}^{p})e_{l}\\
&=&g^{ij}(-J_{k}^{p}\partial_{i}\Gamma_{pj}^{l}+J_{k}^{p}\partial_{p}\Gamma_{ij}^{l}
+J_{p}^{l}\partial_{i}\Gamma_{kj}^{p}-J_{p}^{l}\partial_{k}\Gamma_{ij}^{p})e_{l}+\mathcal{O}(\partial g,\partial J).
\end{eqnarray*}
So we obtain
\begin{displaymath}
(\triangle J+\mathcal{R}+L_{X}J)_{k}^{l}=g^{ij}\partial_{i}\partial_{j}J_{k}^{l}+\mathcal{O}(\partial g,\partial J).
\end{displaymath}
As for the pair $(g^{J},J)$, the lower order terms are still lower order terms, and when we evaluate at $(g_{0},J_{0})$, from the compatibility, we have
\begin{displaymath}
(\mathcal{L}_{(g_{0},0)}\tilde{\mathcal{D}}_{2})(h,K)=\triangle_{g_{0}}K+\mathcal{O}(\partial g,\partial J).
\end{displaymath}
Hence, the total symbol is
\begin{displaymath}
\sigma(\mathcal{L}_{(g_{0},0)}\tilde{\mathcal{D}})(h,K)(x,\xi)=
\left(\begin{array}{cc}
|\xi|^{2}&0\\
0&|\xi|^{2}
\end{array}\right).
\end{displaymath}
By the standard theory of parabolic PDE, there exists a unique short time solution of (\ref{16}).

Next, we show that under (\ref{16}), $(g,J)$ is compatible, where $J=\pi E$. Suppose $(g,J)$ exists for $t\in [0,\epsilon_{0}]$, then by the compactness of $M$, in this time interval, every tensor we involve is bounded. Let $\frac{\partial}{\partial t}J=K$. Then,
\setlength\arraycolsep{2pt}
\begin{eqnarray*}
\frac{\partial}{\partial t}|g^{-J}|_{g^{J}}^{2}&=&2\langle\frac{\partial}{\partial t}(g^{-J}),g^{-J}\rangle_{g^{J}}+C\ast(g^{-J})^{\ast2}\\
&=&2\langle\frac{\partial}{\partial t}\frac{1}{2}(g(\cdot,\cdot)-g(J\cdot,J\cdot)),g^{-J}\rangle_{g^{J}}+C\ast(g^{-J})^{\ast2}\\
&=&2\langle(\frac{\partial}{\partial t}g)^{-J},g^{-J}\rangle_{g^{J}}-\langle g(J\cdot,K\cdot)+g(K\cdot,J\cdot),g^{-J}\rangle_{g^{J}}+C\ast(g^{-J})^{\ast2}\\
&\leq&\langle2(\mathcal{D}_{1}(g^{J},J))^{-J}+2(\triangle_{g_{0}}(g^{-J}))^{-J}-
g(J\cdot,K\cdot)-g(K\cdot,J\cdot),g^{-J}\rangle_{g^{J}}+C|g^{-J}|^{2}_{g^{J}}.
\end{eqnarray*}
Note that $(g^{J},J)$ is compatible and $K=\mathcal{D}_{2}(g^{J},J)$, so by Lemma \ref{compa} and Lemma \ref{04},
\begin{displaymath}
\mathcal{D}_{1}(g^{J},J)^{-J}-\frac{1}{2}(g^{J}(J\cdot,K\cdot)+g^{J}(K\cdot,J\cdot))=0.
\end{displaymath}
So
\setlength\arraycolsep{2pt}
\begin{eqnarray*}
\frac{\partial}{\partial t}|g^{-J}|_{g^{J}}^{2}&\leq&2\langle(\triangle_{g_{0}}(g^{-J}))^{-J}-
g^{-J}(J\cdot,K\cdot)-g^{-J}(K\cdot,J\cdot),g^{-J}\rangle_{g^{J}}+C|g^{-J}|^{2}_{g^{J}}\\
&\leq&2\langle(\triangle_{g_{0}}(g^{-J}))^{-J},g^{-J}\rangle_{g^{J}}+C|g^{-J}|^{2}_{g^{J}}.
\end{eqnarray*}
Noting that $J$ isometrically acts on the space of $\Gamma(T^{*}M\bigotimes^{sym} T^{*}M)$ in the induced metric from $g^{J}$, and $(1,1)$ tensors and $(0,2)+(2,0)$ tensors correspond to $+1$ and $-1$ eigenspace respectively, so they are orthogonal. So
\begin{displaymath}
\langle(\triangle_{g_{0}}(g^{-J}))^{J},g^{-J}\rangle_{g^{J}}=0.
\end{displaymath}
Then,
\begin{displaymath}
\frac{\partial}{\partial t}|g^{-J}|_{g^{J}}^{2}\leq2\langle\triangle_{g_{0}}(g^{-J}),g^{-J}\rangle_{g^{J}}
+C|g^{-J}|^{2}_{g^{J}}.
\end{displaymath}
By definition,
\begin{displaymath}
\triangle_{g_{0}}(g^{-J})=tr_{g_{0}}D_{g_{0}}^{2}(g^{-J}).
\end{displaymath}
Note that the second order term about $g^{-J}$ in $D_{g_{0}}^{2}(g^{-J})$ is the same as in $D_{g^{J}}^{2}(g^{-J})$, so
\begin{displaymath}
\triangle_{g_{0}}(g^{-J})=tr_{g_{0}}(D_{g^{J}}^{2}(g^{-J})+C^{\prime}\ast D_{g^{J}}(g^{-J})+C\ast g^{-J}).
\end{displaymath}
Let $A$ be any tensor, we have the following formula,
\setlength\arraycolsep{2pt}
\begin{eqnarray*}
D^{2}\langle A,A\rangle&=&D(D\langle A,A\rangle)\\
&=&2D(\langle D_{i}A,A\rangle e^{i})\\
&=&2\langle D^{2}_{i,j}A,A\rangle e^{i}\otimes e^{j}+2\langle D_{i}A,D_{j}A\rangle e^{i}\otimes e^{j}.
\end{eqnarray*}
Let $A=g^{-J}$, the metric above be $g^{J}$. And taking trace of each side with respect to $g_{0}$, we obtain
\begin{displaymath}
2\langle tr_{g_{0}}D_{g^{J}}^{2}(g^{-J}), g^{-J}\rangle_{g^{J}}=tr_{g_{0}}D_{g^{J}}^{2}(|g^{-J}|_{g^{J}}^{2})
-2\langle D_{g^{J}}g^{-J}(e_{i}),D_{g^{J}}g^{-J}(e_{j})\rangle_{g^{J}}\langle e^{i},e^{j}\rangle_{g_{0}}.
\end{displaymath}
Note that along this flow, for $t\in [0,\epsilon_{0}]$, $g^{J}$ is uniformly bounded by $g_{0}$, so we have
\begin{displaymath}
2\langle tr_{g_{0}}D_{g^{J}}^{2}(g^{-J}), g^{-J}\rangle_{g^{J}}\leq tr_{g_{0}}D_{g^{J}}^{2}(|g^{-J}|_{g^{J}}^{2})
-2C^{\prime\prime}|D_{g^{J}}g^{-J}|_{g^{J}}^{2}.
\end{displaymath}
Hence,
\begin{displaymath}
\frac{\partial}{\partial t}|g^{-J}|_{g^{J}}^{2}\leq tr_{g_{0}}D_{g^{J}}^{2}(|g^{-J}|_{g^{J}}^{2})
-2C^{\prime\prime}|D_{g^{J}}g^{-J}|_{g^{J}}^{2}+C^{\prime}\ast D_{g^{J}}(g^{-J})\ast g^{-J}+C|g^{-J}|^{2}_{g^{J}}.
\end{displaymath}
By using Cauchy inequality to $C^{\prime}\ast D_{g^{J}}(g^{-J})\ast g^{-J}$, finally we obtain
\begin{displaymath}
\frac{\partial}{\partial t}|g^{-J}|_{g^{J}}^{2}\leq tr_{g_{0}}D_{g^{J}}^{2}(|g^{-J}|_{g^{J}}^{2})+C|g^{-J}|^{2}_{g^{J}}.
\end{displaymath}
Notice that $tr_{g_{0}}D_{g^{J}}^{2}$ is elliptic and $|g^{-J}|^{2}=0$ at $t=0$.
Then by maximal principle, considering $e^{-Ct}|g^{-J}|^{2}$, we have $|g^{-J}|^{2}=0$ for $t\in [0,\epsilon_{0}]$, i.e. $(g,J)$ is compatible.
Since $\epsilon_{0}$ is arbitrary, $(g,J)$ is always compatible as long as the solution exists.
Because the positivity of $g$ is a open condition, we may assume $g$ is positive in short time.
Then the short time solution of (\ref{16}) gives the short time solution of (\ref{17}).

Now, let $(\tilde{g}(t),\tilde{J}(t))$ be a solution of (\ref{17}), $\varphi_{t}$ be the one-parameter family of diffeomorphisms generated by $-X(t)$ defined as above. Let $g(t)=\varphi_{t}^{\ast}\tilde{g}(t)$, $J(t)=\varphi_{t}^{\ast}\tilde{J}(t)$, then,
\setlength\arraycolsep{2pt}
\begin{eqnarray}
\frac{\partial}{\partial t}g&=&\frac{\partial}{\partial t}(\varphi_{t}^{\ast}\tilde{g}(t))\nonumber\\
&=&\varphi_{t}^{\ast}(\frac{\partial}{\partial t}\tilde{g}(t)+L_{(-X(t))}\tilde{g}(t))\nonumber\\
&=&\varphi_{t}^{\ast}(-2Ric(\tilde{g}(t))
+Q_{1}(\tilde{g}(t)))\label{phi}\\
&=&-2Ric(\varphi_{t}^{\ast}\tilde{g}(t))
+Q_{1}(\varphi_{t}^{\ast}\tilde{g}(t))\nonumber\\
&=&-2Ric(g)+Q_{1}(g)\nonumber.
\end{eqnarray}
So $g(t)$ satisfies the equation. Similarly, $J(t)$ also satisfies the equation. And $(g(t),J(t))$ differs from $(\tilde{g}(t),\tilde{J}(t))$ by a diffeomorphism, so $(g(t),J(t))$ is also an almost Hermitian pair. So we finish the existence part of the theorem.

For uniqueness, let $(g_{i},J_{i})$ be two solutions of (\ref{flow}), $i=1,2$. Since $M$ is compact, we can solve the harmonic heat flow $\phi_{i}(t)$ for short time,
\begin{eqnarray*}
\frac{\partial}{\partial t}\phi_{i}(t)&=&\triangle_{g_{i},\overline{g}}\phi_{i}(t)\\
\phi_{i}(0)&=&id,
\end{eqnarray*}
where $\overline{g}$ is the same fixed metric as above. We can also assume $\phi_{i}(t)$ are diffeomorphisms. Let $\hat{g}_{i}=(\phi_{i}^{-1}(t))^{\ast}g_{i}(t)$. Note that
\setlength\arraycolsep{2pt}
\begin{eqnarray*}
(\frac{\partial}{\partial t}\phi_{i})(p)&=&(\triangle_{g_{i},\overline{g}}\phi_{i})(p)\\
&=&(\triangle_{\hat{g}_{i},\overline{g}}id)(\phi_{i}(p))\\
&=&(-\hat{g}^{ij}(\hat{\Gamma}_{ij}^{k}-\overline{\Gamma}_{ij}^{k})\frac{\partial}{\partial x^{k}})(\phi_{i}(p))\\
&=&-X_{\hat{g}}(\phi_{i}(p)).
\end{eqnarray*}
Then taking time derivative of $(\phi_{i}(t))^{\ast}\hat{g}_{i}(t)=g_{i}(t)$, and doing the similar calculation in (\ref{phi}),
we see that both $\hat{g}_{i}(t)$ satisfy (\ref{17}) and they share the same initial data.
Since we have proved the compatibility, the symbol of (\ref{17}) is $Id$ as we calculated, so the solution of (\ref{17}) is unique. Then, we obtain
\begin{displaymath}
\hat{g}_{1}(t)=\hat{g}_{2}(t)=\hat{g}(t),\quad \hat{J}_{1}(t)=\hat{J}_{2}(t)=\hat{J}(t).
\end{displaymath}
Then from the uniqueness of
\setlength\arraycolsep{2pt}
\begin{eqnarray*}
\frac{\partial}{\partial t}\phi(t)&=&-X_{\hat{g}}(\phi(t))\\
\phi(0)&=&id,
\end{eqnarray*}
we see the uniqueness of $(g,J)$ for a short while. Then by continuity, $(g,J)$ is unique as long as it exists.

Next, we check two special cases. Suppose the initial data is almost K\"{a}hler, then we run the symplectic curvature flow (\ref{symp}).
%\setlength\arraycolsep{2pt}
%\begin{eqnarray*}
%\frac{\partial}{\partial t}g&=&-2Ric+\frac{1}{2}B^{1}-B^{2}\\
%\frac{\partial}{\partial t}J&=&\triangle{J}+\mathcal{N}+\mathcal{R}\\
%g(0)&=&g_{0}\\
%J(0)&=&J_{0}.
%\end{eqnarray*}
By definitions and Lemma \ref{s},
we see in this situation, $(g,J)$ also satisfies (\ref{flow}). So from the uniqueness of (\ref{flow}),
if the initial data is almost K\"{a}hler, then (\ref{flow}) coincides with symplectic curvature flow.
And the similar argument holds for pluriclosed case when we apply Lemma \ref{p}.

%Similarly, if the initial data is pluriclosed, then we run the equivalent pluriclosed flow (\ref{pluri}).
%\setlength\arraycolsep{2pt}
%\begin{eqnarray*}
%\frac{\partial}{\partial t}g&=&-2Ric+\frac{1}{2}H^{2}\\
%\frac{\partial}{\partial t}J&=&\triangle{J}+\mathcal{R}+Q\\
%g(0)&=&g_{0}\\
%J(0)&=&J_{0}.
%\end{eqnarray*}
%By definitions and Lemma \ref{35}, we see that in this situation, $(g,J)$ also satisfies (\ref{flow}).
%So from the uniqueness of (\ref{flow}), if the initial data is pluriclosed, then (\ref{flow}) is equivalent to pluriclosed flow.
Finally, we prove flow (\ref{flow}) preserves the integrability of $J$. Let $(g_{0},J_{0})$ be an Hermitian structure. Fix $J_{0}$, consider the following flow,
\begin{eqnarray*}
\frac{\partial}{\partial t}\tilde{g}&=&-2Ric_{\tilde{g}}+Q_{1}(\tilde{g},J_{0})-L_{\theta^{\sharp}(\tilde{g},J_{0})}\tilde{g}\\
\tilde{g}(0)&=&g_{0}.
\end{eqnarray*}
By DeTurck trick, we see $\tilde{g}(t)$ exists for a while, but not necessary compatible with $J_{0}$ now.
Then by a gauge transformation induced by $\theta^{\sharp}(\tilde{g},J_{0})$, we obtain a short time solution $(g(t),J(t))$ for the following flow,
\begin{eqnarray*}
\frac{\partial}{\partial t}g&=&-2Ric_{g}+Q_{1}(g,J)\\
\frac{\partial}{\partial t}J&=&L_{\theta^{\sharp}(g,J)}J\\
g(0)&=&g_{0}\\
J(0)&=&J_{0}.
\end{eqnarray*}
Notice that we still don't know the compatibility of $(g,J)$ now, but since $J$ is changed just by a diffeomorphism, $N$ always vanishes. By Lemma \ref{gau}, one may write $Q_{2}-\mathcal{Q}+\mathcal{N}$ in terms of $N$ in the almost Hermitian setting. We denote such a tensor $N_{0}$, i.e. $N_{0}$ is in terms of $N$ and when $(g,J)$ is compatible, $N_{0}=Q_{2}-\mathcal{Q}+\mathcal{N}$ . So the above flow is the same as following,
\begin{eqnarray*}
\frac{\partial}{\partial t}g&=&-2Ric_{g}+Q_{1}(g,J)\\
\frac{\partial}{\partial t}J&=&L_{\theta^{\sharp}(g,J)}J+N_{0}(g,J)-\overline{N}(g,J)-\mathcal{K}(g,J)\\
g(0)&=&g_{0}\\
J(0)&=&J_{0}.
\end{eqnarray*}
Then by Lemma \ref{Lee}, and using the same argument in the proof of short time existence above, one see $(g,J)$ is compatible and coincides with (\ref{flow}), so the integrability of $J$ is preserved.

Therefore we finish the proof of Theorem \ref{main}.
\qed
\begin{rem}
In \cite{ST symp}, Streets and Tian introduced almost Hermitian curvature flow, where the symbol term they deforming $J$ is $-\mathcal{K}$. From Lemma \ref{Lee}, we see modulo lower order terms, $-\mathcal{K}$ differs from $\triangle J+\mathcal{R}$ just by a gauge term. While if we also change the evolution of $g$ by the same gauge transformation, the second derivative of $g$ will appear in $L_{\theta^{\sharp}}g$. So in general, our flow is not in the family of almost Hermitian curvature flow.
\end{rem}

\section{Proof of Theorem \ref{main2} and Theorem \ref{main3}}\label{pf23}
First, we derive the evolution equations of $DJ$, $Rm$ and their higher covariant derivatives.
\begin{lem}\label{edj}
Under (\ref{flow}),
\begin{eqnarray*}
\frac{\partial}{\partial t}DJ=\triangle DJ+Rm*DJ+J^{*2}*DJ^{*3}+J^{*3}*DJ*D^{2}J.
\end{eqnarray*}
\end{lem}
\pf
Using the fact $\triangle DT-D\triangle T=DRm*T+Rm*DT$, we have
\begin{eqnarray*}
\frac{\partial}{\partial t}DJ&=&\dot{\Gamma}*J+D\dot{J}\\
&=&D(Rm+J^{*2}*DJ^{*2})*J+D(\triangle J+Rm*J+J*DJ^{*2})\\
&=&\triangle DJ+DRm*J+Rm*DJ+J^{*2}*DJ^{*3}+J^{*3}*DJ*D^{2}J.
\end{eqnarray*}
Hence we only need to show there is no $DRm*J$ term. It is the same calculation as in \cite{ST symp},
since the only differences are the first order terms of $J$, which does not involve $DRm$ term.
\qed
\begin{lem}
Under (\ref{flow}),
\begin{eqnarray*}
\frac{\partial}{\partial t}Rm=\triangle Rm+Rm^{*2}+Rm*J^{*2}*DJ^{*2}+\sum_{
\tiny{
\begin{array}{l}
0\leq k_{1},\ldots, k_{4}\leq 3\\
k_{1}+\cdots+k_{4}=4
\end{array}}
}D^{k_{1}}J*\cdots*D^{k_{4}}J.
\end{eqnarray*}
\end{lem}
\pf
From the variation formula in Ricci flow, see \cite{CK}, we have, let $\frac{\partial}{\partial t}g=h$,
\begin{eqnarray*}
\frac{\partial}{\partial t}Rm(X,Y,Z,W)&=&\frac{1}{2}(h(Rm(X,Y)Z,W)-h(Rm(X,Y)W,Z))\\
&&+\frac{1}{2}(D^{2}_{Y,W}h(X,Z)-D^{2}_{X,W}h(Y,Z)+D^{2}_{X,Z}h(Y,W)-D^{2}_{Y,Z}h(X,W)).
\end{eqnarray*}
And when $h=-2Ric$,
\begin{eqnarray*}
\frac{\partial}{\partial t}Rm=\triangle Rm+Rm^{*2}.
\end{eqnarray*}
Notice that in (\ref{flow}), $h=\frac{\partial}{\partial t}g=-2Ric+J^{*2}*DJ^{*2}$, then we obtain the evolution equation of $Rm$.
\qed
\begin{prop}\label{evo}
Under (\ref{flow}),
\begin{eqnarray*}
\frac{\partial}{\partial t}D^{k}J&=&\triangle D^{k}J+\sum_{
\tiny{
\begin{array}{l}
l_{1}+\cdots+l_{5}=k+2\\
0\leq l_{1}, \ldots, l_{5}\leq k+1
\end{array}}}
D^{l_{1}}J*\cdots*D^{l_{5}}J+\sum_{l=0}^{k-1}{D^{l}Rm*D^{k-l}}J,\\
\frac{\partial}{\partial t}D^{k}Rm&=&\triangle D^{k}Rm+\sum_{
\tiny{
\begin{array}{l}
l_{1}+\cdots+l_{4}=k+4\\
0\leq l_{1}, \ldots, l_{4}\leq k+3
\end{array}}}
D^{l_{1}}J*\cdots*D^{l_{4}}J
+\sum_{l=0}^{k}
D^{l}Rm*D^{k-l}Rm\\
&&+\sum_{
\begin{array}{l}
0\leq l_{0}\leq k
\end{array}}
\sum_{
\tiny{
\begin{array}{l}
l_{1}+\cdots+l_{4}=k+2-l_{0}\\
0\leq l_{1},\ldots,l_{4}\leq k+1
\end{array}}}
D^{l_{0}}Rm*D^{l_{1}}J*\cdots*D^{l_{4}}J.
\end{eqnarray*}
\end{prop}
\pf
By using Lemma \ref{edj} and the fact $\frac{\partial}{\partial t}\Gamma=D(Rm+J^{*2}*DJ^{*2})$, we have
\begin{eqnarray*}
\frac{\partial}{\partial t}D^{k}J&=&\frac{\partial}{\partial t}\Gamma*D^{k-1}J+D\frac{\partial}{\partial t}D^{k-1}J\\
&=&\sum_{l=0}^{k-2}D^{l}\frac{\partial}{\partial t}\Gamma*D^{k-1-l}J+D^{k-1}\frac{\partial}{\partial t}DJ\\
&=&\sum_{l=0}^{k-2}D^{l}D(Rm+J^{*2}*DJ^{*2})*D^{k-1-l}J\\
&&+D^{k-1}(\triangle DJ+Rm*DJ+J^{*2}*DJ^{*3}+J^{*3}*DJ*D^{2}J).
\end{eqnarray*}
Interchanging $D$ and $\triangle$, then we observe that the highest order of $Rm$ is $k-1$,
and the highest order of $J$ is $k+1$ if not involving $Rm$. Then we obtain the evolution equation of $D^{k}J$.

As for the evolution equation of $D^{k}Rm$, the calculation is similar. The key point is to observe the highest order.
\qed

Now we can use Proposition \ref{evo} to prove Theorem \ref{main2} and Theorem \ref{main3}.\\
\textit{Proof of Theorem \ref{main2}. }
The proof is similar to the higher derivative estimates in Ricci flow \cite{CK}.

We assume $t|D^{2}J|\leq C$ first. By induction, we will prove
\begin{eqnarray*}
(P):\quad|D^{k}J|\leq\frac{C}{t^{\frac{k}{2}}}, |D^{k-2}Rm|\leq\frac{C}{t^{\frac{k}{2}}}.
\end{eqnarray*}
$(P)$ holds when $k=2$ from the assumption.

Now we assume $(P)$ holds for $k-1$.
Consider
\begin{eqnarray*}
F(t)=t^{k+1}(|D^{k}J|^{2}+|D^{k-2}Rm|^{2})+\lambda t^{k}(|D^{k-1}J|^{2}+|D^{k-3}Rm|^{2}),
\end{eqnarray*}
where $\lambda$ is a large constant to be determined. We will show
\begin{eqnarray}
\frac{\partial}{\partial t}F\leq\triangle F+C.\label{est}
\end{eqnarray}
Then by maximum principle, $(P)$ holds for $k$. Now we prove (\ref{est}), by using Proposition \ref{evo},
\begin{eqnarray*}
\frac{\partial}{\partial t}|D^{k}J|^{2}&=&(Rm+J^{*2}*DJ^{*2})*D^{k}J^{*2}\\
&&+2\langle D^{k}J, \quad\triangle D^{k}J+\sum_{
\tiny{
\begin{array}{l}
l_{1}+\cdots+l_{5}=k+2\\
0\leq l_{1}, \ldots, l_{5}\leq k+1
\end{array}}}
D^{l_{1}}J*\cdots*D^{l_{5}}J+\sum_{l=0}^{k-1}{D^{l}Rm*D^{k-l}}J\rangle\\
&=&(Rm+J^{*2}*DJ^{*2})*D^{k}J^{*2}+\triangle|D^{k}J|^{2}-2|D^{k+1}J|^{2}\\
&&+D^{k}J*(\sum_{
\tiny{
\begin{array}{l}
l_{1}+\cdots+l_{5}=k+2\\
0\leq l_{1}, \ldots, l_{5}\leq k+1
\end{array}}}
D^{l_{1}}J*\cdots*D^{l_{5}}J+\sum_{l=0}^{k-1}{D^{l}Rm*D^{k-l}}J)\\
&=&\triangle|D^{k}J|^{2}-2|D^{k+1}J|^{2}+(Rm+J^{*2}*DJ^{*2})*D^{k}J^{*2}\\
&&+D^{k}J*D^{k+1}J*DJ*J^{*3}+D^{k}J*D^{k}J*DJ^{*2}*J^{*2}+D^{k}J*D^{k}J*D^{2}J*J^{*3}\\
&&+D^{k}J*\sum_{
\tiny{
\begin{array}{l}
l_{1}+\cdots+l_{5}=k+2\\
0\leq l_{1}, \ldots, l_{5}\leq k-1
\end{array}}}
D^{l_{1}}J*\cdots*D^{l_{5}}J\\
&&+D^{k}J*Rm*D^{k}J+D^{k}J*D^{k-1}Rm*DJ+D^{k}J*D^{k-2}Rm*D^{2}J\\
&&+D^{k}J*\sum_{l=1}^{k-3}{D^{l}Rm*D^{k-l}}J.
\end{eqnarray*}
From the assumption,
\begin{eqnarray*}
\frac{\partial}{\partial t}|D^{k}J|^{2}&\leq&\triangle|D^{k}J|^{2}-2|D^{k+1}J|^{2}+\frac{C}{t}|D^{k}J|^{2}+\frac{C}{t^{\frac{1}{2}}}|D^{k}J||D^{k+1}J|\\
&&+\frac{C}{t^{\frac{k+2}{2}}}|D^{k}J|+\frac{C}{t^{\frac{1}{2}}}|D^{k}J||D^{k-1}Rm|+\frac{C}{t}|D^{k}J||D^{k-2}Rm|.
\end{eqnarray*}
Similarly, we obtain
\begin{eqnarray*}
\frac{\partial}{\partial t}|D^{k-2}Rm|^{2}&\leq&\triangle|D^{k-2}Rm|^{2}-2|D^{k-1}Rm|^{2}+\frac{C}{t}|D^{k-2}Rm|^{2}+\frac{C}{t^{\frac{1}{2}}}|D^{k-2}Rm||D^{k+1}J|\\
&&+\frac{C}{t^{\frac{k+2}{2}}}|D^{k-2}Rm|+\frac{C}{t}|D^{k}J||D^{k-2}Rm|.
\end{eqnarray*}
Then by Cauchy-Schwarz inequality,
\begin{eqnarray*}
\frac{\partial}{\partial t}(t^{k+1}(|D^{k}J|^{2}+|D^{k-2}Rm|^{2}))&\leq&\triangle(t^{k+1}(|D^{k}J|^{2}+|D^{k-2}Rm|^{2}))
-t^{k+1}(|D^{k+1}J|^{2}+|D^{k-1}Rm|^{2})\\
&&+Ct^{k}(|D^{k}J|^{2}+|D^{k-2}Rm|^{2})+C.
\end{eqnarray*}
Replacing $k$ by $k-1$ and using assumption, we obtain
\begin{eqnarray*}
\frac{\partial}{\partial t}(t^{k}(|D^{k-1}J|^{2}+|D^{k-3}Rm|^{2}))&\leq&\triangle(t^{k}(|D^{k-1}J|^{2}+|D^{k-3}Rm|^{2}))
-t^{k}(|D^{k}J|^{2}+|D^{k-2}Rm|^{2})\\
&&+Ct^{k-1}(|D^{k-1}J|^{2}+|D^{k-3}Rm|^{2})+C\\
&\leq&\triangle(t^{k}(|D^{k-1}J|^{2}+|D^{k-3}Rm|^{2}))-t^{k}(|D^{k}J|^{2}+|D^{k-2}Rm|^{2})+C.
\end{eqnarray*}
Then
\begin{eqnarray*}
\frac{\partial F}{\partial t}&\leq& \triangle F-t^{k+1}(|D^{k+1}J|^{2}+|D^{k-1}Rm|^{2})+(C-\lambda)t^{k}(|D^{k}J|^{2}+|D^{k-2}Rm|^{2})+C\\
&\leq&\triangle F+(C-\lambda)t^{k}(|D^{k}J|^{2}+|D^{k-2}Rm|^{2})+C.
\end{eqnarray*}
We choose $\lambda=C$, then (\ref{est}) holds.

Now, we prove $t|D^{2}J|\leq C$. For $p\in M$, if $|D^{2}J|_{p,t}\neq 0$, then similarly, by Proposition \ref{evo},
\begin{eqnarray*}
\frac{\partial}{\partial t}|D^{2}J|&=&\frac{1}{2|D^{2}J|}\frac{\partial}{\partial t}|D^{2}J|^{2}\\
&=&\frac{1}{2|D^{2}J|}(\triangle|D^{2}J|^{2}-2|D^{3}J|^{2}+D^{2}J^{*3}*J^{*3}+D^{3}J*D^{2}J*DJ*J^{*3}\\
&&+D^{2}J^{*2}*DJ^{*2}*J^{*2}+D^{2}J*DJ^{*4}+D^{2}J^{*2}*Rm+D^{2}J*DJ*DRm).
\end{eqnarray*}
Notice that for $|D^{2}J|_{p,t}\neq 0$,
\begin{eqnarray*}
\triangle|D^{2}J|^{2}=2|D^{2}J|\triangle|D^{2}J|+2|D|D^{2}J||^{2}.
\end{eqnarray*}
So,
\begin{eqnarray*}
\frac{\partial}{\partial t}|D^{2}J|&=&\triangle|D^{2}J|+\frac{|D|D^{2}J||^{2}}{|D^{2}J|}+\frac{1}{2|D^{2}J|}(-2|D^{3}J|^{2}+D^{2}J^{*3}*J^{*3}+D^{3}J*D^{2}J*DJ*J^{*3}\\
&&+D^{2}J^{*2}*DJ^{*2}*J^{*2}+D^{2}J*DJ^{*4}+D^{2}J^{*2}*Rm+D^{2}J*DJ*DRm)\\
&\leq&\triangle|D^{2}J|+\frac{|D|D^{2}J||^{2}}{|D^{2}J|}-\frac{|D^{3}J|^{2}}{|D^{2}J|}+C(|D^{2}J|^{2}+\frac{|D^{3}J|}{t^{\frac{1}{2}}}+\frac{|D^{2}J|}{t}+\frac{1}{t^{2}}+\frac{|DRm|}{t^{\frac{1}{2}}}).
\end{eqnarray*}
Consider
\begin{eqnarray*}
G(t)=t^{2}|D^{2}J|+\mu t^{2}|DJ|^{2}+t^{3}|Rm|^{2},
\end{eqnarray*}
where $\mu$ is a large constant to be determined.

Then for $|D^{2}J|\neq 0$,
\begin{eqnarray*}
\frac{\partial}{\partial t}G&\leq&\triangle G-t^{2}\frac{|D^{3}J|^{2}}{|D^{2}J|}-2\mu t^{2}|D^{2}J|^{2}-2t^{3}|DRm|^{2}\\
&&+C(t^{2}|D^{2}J|^{2}+t^{\frac{3}{2}}|D^{3}J|+\mu t|D^{2}J|+\mu+t^{\frac{3}{2}}|DRm|)+\langle D|t^{2}D^{2}J|, \frac{D|D^{2}J|}{|D^{2}J|}\rangle\\
&\leq&\triangle G-\frac{1}{2}t^{2}\frac{|D^{3}J|^{2}}{|D^{2}J|}-\frac{1}{2}t^{2}|D^{2}J|^{2}-\frac{1}{2}t^{3}|DRm|^{2}
+\langle D|t^{2}D^{2}J|, \frac{D|D^{2}J|}{|D^{2}J|}\rangle+C,
\end{eqnarray*}
where $\mu$ is determined now.

Then
\begin{eqnarray*}
\frac{\partial}{\partial t}G&\leq&\triangle G-\frac{1}{2}t^{2}\frac{|D^{3}J|^{2}}{|D^{2}J|}-\frac{1}{2} t^{2}|D^{2}J|^{2}-\frac{1}{2}t^{3}|DRm|^{2}+C\\
&&+\langle DG, \frac{D|D^{2}J|}{|D^{2}J|}\rangle-\mu t^{2}\langle D|DJ|^{2}, \frac{D|D^{2}J|}{|D^{2}J|}\rangle-t^{3}\langle D|Rm|^{2}, \frac{D|D^{2}J|}{|D^{2}J|}\rangle.
\end{eqnarray*}
Notice that
\begin{eqnarray*}
|D|DJ|^{2}|&\leq&|2\langle DDJ, DJ\rangle|\leq2|D^{2}J||DJ|,\\
|D|D^{2}J||&=&\frac{|D|D^{2}J|^{2}|}{2|D^{2}J|}\leq|D^{3}J|.
\end{eqnarray*}
Hence,
\begin{eqnarray*}
\frac{\partial}{\partial t}G&\leq&\triangle G-\frac{1}{4}t^{2}\frac{|D^{3}J|^{2}}{|D^{2}J|}-\frac{1}{4} t^{2}|D^{2}J|^{2}-\frac{1}{2}t^{3}|DRm|^{2}+C+\langle DG, \frac{D|D^{2}J|}{|D^{2}J|}\rangle+C\frac{t^{2}|DRm|^{2}}{|D^{2}J|}.
\end{eqnarray*}
So if we suppose $|D^{2}J|\geq \frac{4C}{t}$, we have the estimate,
\begin{eqnarray*}
\frac{\partial}{\partial t}G&\leq&\triangle G+\langle DG, \frac{D|D^{2}J|}{|D^{2}J|}\rangle+C,
\end{eqnarray*}
where $C=C(n,K)$. That is to say, for any $(p,t)$, either we have the estimate $|D^{2}J|\leq \frac{4C}{t},$
\begin{eqnarray*}
\text{or} \quad \frac{\partial}{\partial t}G\leq\triangle G+\langle DG, \frac{D|D^{2}J|}{|D^{2}J|}\rangle+C.
\end{eqnarray*}
Let $\overline{G}=G-Ct$, where $C$ is chosen suitable. We obtain that either $\overline{G}\leq 0$,
\begin{eqnarray*}
\text{or} \quad \frac{\partial}{\partial t}\overline{G}\leq\triangle \overline{G}+\langle D\overline{G}, \frac{D|D^{2}J|}{|D^{2}J|}\rangle.
\end{eqnarray*}
Notice that $\overline{G}=0$ when $t=0$. Then one may apply maximum principle to show that $\overline{G}\leq 0$ for every $(p,t)$, which implies the desired estimate. So we finish the proof.
\qed
\begin{rem}
Theorem \ref{main2} is scaling invariant when we replace $g(t)$ by $\overline g(t)=cg(\frac{t}{c})$.
\end{rem}
\textit{Proof of Theorem \ref{main3}. }
The argument is standard as in Ricci flow \cite{CK}. We just sketch the proof.

Suppose not, then $|Rm|, |DJ|$ are bounded. From Theorem \ref{main2}, all covariant derivatives of $Rm$ and $J$ are bounded. Then we see $g$ are uniformly bounded. We fix a coordinate atlas. From the evolution equation of $\Gamma$ and the boundedness of covariant derivatives of $Rm$ and $J$, we obtain the boundedness of $\Gamma$.
Then we obtain the boundedness of $\partial g, \partial J$ and by induction we see that $\partial^{k} g, \partial^{k} J$ and $\partial^{k}\Gamma$ are bounded.
Finally, we obtain that $\frac{\partial^{l}}{\partial t^{l}}\partial^{k} g, \frac{\partial^{l}}{\partial t^{l}}\partial^{k} J$ are bounded.
Then by theorems in mathematical analysis, $(g(t), J(t))$ can be extended to $(g(T), J(T))$ smoothly in all variables of space and time.
The almost Hermitian condition is guaranteed by the continuity. Then from the short time existence, $(g(t), J(t))$ exists for $t\in[0,T+\epsilon)$,
which is a contradiction to the maximality of $T$.
\qed


\begin{thebibliography}{99}
\bibitem{B} Boling J. $Homogeneous~ solutions~ of~ pluriclosed~ flow~ on~ closed~ complex~ surfaces$, arXiv preprint arXiv:1404.7106, 2014.
\bibitem{CK} Chow B., Knopf, D. $The~ Ricci~ flow:~ an~ introduction$, AMS Bookstore, 2004.
\bibitem{EFV} Enrietti, N., Fino, A., Vezzoni, L. $The~ pluriclosed~ flow~ on~ nilmanifolds~ and~ Tamed~ symplectic~ forms$, arXiv preprint arXiv:1210.4816, 2012.
\bibitem{E} Enrietti, N. $Static~ SKT~ metrics~ on~ Lie~ groups$, Manuscripta Mathematica, 2013: 1-15.
\bibitem{F} Fern\'{a}ndez-Culma, E. $Soliton~ almost~ K\ddot{a}hler~ structures~ on~ 6-dimensional~ nilmanifolds~ for~ the~
symplectic~ \\curvature~ flow$, arXiv preprint arXiv:1303.5461, 2013.
\bibitem{Ga} Gauduchon, P. $Hermitian~ connections~ and~ Dirac~ operators$, Bollettino della Unione Matematica Italiana-B, 1997 (2): 257-288.
\bibitem{G} Gualtieri, M. $Generalized~ complex~ geometry$, arXiv preprint math/0401221, 2004.
\bibitem{H} Hitchin, N. $Generalized~ Calabi-Yau~ manifolds$, The Quarterly Journal of Mathematics, 2003, 54(3): 281-308.
%\bibitem{KN} Kobayashi, S., Nomizu, K. $Foundations of differential geometry Volume I[M]$, Wiley-Interscience, 1996.
\bibitem{LW} Lauret, J., Will, C. $On~ the~ symplectic~ curvature~ flow~ for~ locally~ homogeneous~ manifolds$ arXiv preprint arXiv:1405.6065, 2014.
\bibitem{NN} Newlander, A., Nirenberg, L. $Complex~ analytic~ coordinates~ in~ almost~ complex~ manifolds$, The Annals of Mathematics, 1957, 65(3): 391-404.
\bibitem{Pe} Perelman, G. $The~ entropy~ formula~ for~ the~ Ricci~ flow~ and~ its~ geometric~ applications$, arXiv preprint math/0211159, 2002.
\bibitem{Po} Pook, J. $Homogeneous~ and~ locally~ homogeneous~ solutions~ to~ symplectic~ curvature~ flow$, arXiv preprint arXiv:1202.1427, 2012.
\bibitem{Sm} Smith, D. $Stability~ of~ the~ Almost~ Hermitian~ Curvature~ Flow$, arXiv preprint arXiv:1308.6214, 2013.
\bibitem{Streets pgs} Streets, J. $Pluriclosed~ flow~ on~ generalized~ K\ddot{a}hler~ manifolds~ with~ split~ tangent~ bundle$, arXiv preprint arXiv:1405.0727, 2014.
\bibitem{Streets pbg} Streets, J. $Pluriclosed~ flow, ~Born-Infeld ~geometry,~ and~ rigidity~ results~ for ~ generalized~ K\ddot{a}hler~ manifolds$, arXiv preprint arXiv:1502.02584, 2015.
\bibitem{ST pluri} Streets, J., Tian, G. $A~ parabolic~ flow~ of~ pluriclosed~ metrics$, Int. Math. Res. Notices (2010), (16): 3101-3133.
\bibitem{ST reg pluri} Streets, J., Tian, G. $Regularity~ results~ for~ pluriclosed~ flow$, to appear Geometry \& Topology.
\bibitem{ST her} Streets, J., Tian, G. $Hermitian~ curvature~ flow$, Journal of the European Mathematical Society, 2011, 13(3): 601-634.
\bibitem{ST symp} Streets, J., Tian, G. $Symplectic~ curvature~ flow$, Journal f\"{u}r die reine und angewandte Mathematik (Crelles Journal), 2011.
\bibitem{ST gkg} Streets, J., Tian, G. $Generalized~ K\ddot{a}hler~ geometry~ and~ the~ pluriclosed~ flow$, Nuclear Physics B, 2012, 858(2): 366-376.
\bibitem{SW} Streets, J., Warren, M. $Evans-Krylov~ Estimates~ for~ a~ nonconvex~ Monge-Amp\grave{e}re~ equation$, arXiv preprint arXiv:1410.2911, 2014.
\bibitem{V} Vezzoni, L. $On~ Hermitian~ curvature~ flow~ on~ almost~ complex~ manifolds$, Differential Geometry and its Applications, 2011, 29(5): 709-722.
\end{thebibliography}
\end{document}